\documentclass[12pt]{amsart}

\usepackage{graphicx}
\usepackage{color}
\usepackage{amssymb}
\usepackage{latexsym}

\newcommand{\R}{\mathbb{R}}
\newcommand{\C}{\mathbb{C}}
\newcommand{\N}{\mathbb{N}}
\newcommand{\T}{\mathbb{T}}
\newcommand{\Z}{\mathbb{Z}}
\newcommand{\be}{\begin{equation}}
\newcommand{\ee}{\end{equation}}
\def\into{\int \hspace*{-3mm} - \,}

\newtheorem{thm}{Theorem}[section]
\newtheorem{lem}[thm]{Lemma}
\newtheorem{coro}[thm]{Corollary}

\newtheorem{rem}[thm]{Remark}

\author[Luc Molinet]{Luc  Molinet}

\title[Ill-posedness for periodic cubic Schr\"odinger equation]{On ill-posedness for
the one-dimensional periodic  cubic Schrodinger equation}
\address[Luc Molinet] {LAGA, Institut Galil\'ee, Universit\'e Paris 13 \\
93430 Villetaneuse}

\begin{document}

\begin{abstract}
We prove the ill-posedness  in $ H^s(\T) $,
$s<0$, of the periodic cubic Schr\"odinger equation in the sense that the flow-map is not continuous 
 from $H^s(\T) $ into itself for any fixed $ t\neq 0 $. This result
  is  slightly stronger  than  the one in  \cite{CCT2} where the discontinuity of the solution map 
    is established.  Moreover  our proof is   different and
  clarifies the ill-posedness phenomena. Our approach relies on a new result on the behavior of the associated
   flow-map with respect to the weak topology of $ L^2(\T) $.

\end{abstract}

\subjclass[2000]{}

\date{}

\maketitle

\section{Introduction}
In this work we are interested in the failure of well-posedness in some spaces of rough functions
 for the one dimensional cubic Schr\"odinger  equation in the periodic setting.
 Throughout this paper we will say that  a Cauchy problem is  (locally)  well-posed   in some  normed
  function space $ X $  if, for any initial data $ u_0\in X $, there exist a radius $ R>0 $, a  time
  $ T>0 $  and a unique solution $ u$ ,  belonging to some space-time function space
  continuously embedded in $ C([0,T];X) $, such that for any $ t\in [0,T] $ the map 
   $ u_0\mapsto u(t) $ is continuous from the ball of $ X $ centered at $ u_0 $ with radius $ R $
    into $ X $.    A Cauchy problem will be  said to be ill-posed if it is not well-posed.

 Let us recall that the cubic Schr\"odinger equation reads
   \begin{equation}
 i u_t +u_{xx}+\gamma |u|^2 u=0 \label{schro3} ,
 \end{equation}
 where $ \gamma \neq 0 $ is a real number. In the one-dimensional periodic setting,
 $ u=u(t,x) $ is a function from $ \R\times \T$   to $\C $ where $\T:=\R/2\pi \Z $. Note that this equation is completely integrable (cf. \cite{A}) even if we will not  directly exploit this particularity.

  It is well-known that this equation
   enjoys two symmetries. Namely, the dilation symmetry : $u(t,x)\mapsto
    \lambda u(\lambda^2 t, \lambda x ) $ and the Galilean invariance : $ u(t,x)\mapsto
     e^{i\alpha x/2 } e^{-i\alpha^2 t/4} u(t,x-\alpha t) $.
Recall that  to each equation that enjoys a dilation symmetry, one can  associate  a so called critical scaling Sobolev index $s_c $. This is the index of the homogeneous Sobolev space on the line whose semi-norm remains
 invariant by
  the dilation symmetry of the equation. It is commonly conjectured (at least when working on $\R $ or $ \T$)
   that an equation must be ill-posed in $ H^s $ as soon as $ s<s_c$ (However, as far as the author knows,
    there exists no general proof of this fact). On the other hand, one can expect that the equations
     are well-posed above $s_c $. This is indeed the case for a large class of parabolic equations as the
     dissipative Burgers equation, on $ \R $ or $\T $, and for dispersive equations  as the nonlinear
      Schr\"odinger equation on $ \R^N $ and the generalized Korteweg-de Vries equation on $ \R $
 with $L^2 $-critical or $ L^2 $-supercritical  nonlinearity.

For the cubic Schr\"odinger equation $ s_c=-1/2 $ and
well-posedness is known to hold in $ H^s(\R) $ and $ H^s(\T) $
 for $ s\ge 0 $ (see \cite{Ts} on $ \R $ and \cite{Bo1} on $ \T $). However, as we mentioned, there
  exists another symmetry, the Galilean invariance, for which the critical index is $ 0 $. This suggests
  that this equation may be ill-posed in $ H^s $ as soon as $ s<0 $. A first step in this direction was due
  to  Bourgain (\cite{Bo4}) who observed that  for a wide class of dispersive equations  there exists
   a critical index $ s^\infty_c >s_c $ below which the flow-map of the equation (if it exists) cannot
    be smooth  whereas it is real analytic\footnote{For dispersive periodic equations   whose nonlinear term
     is of the form $ u^q u_x $ the smoothness of the flow-map holds  not for the original equation but
     for the equation satisfied by $ \tilde{u}(t,x)=u(t,x-\int_0^t \into u^q) $. Therefore the non-smoothness
     of the flow-map has to be shown on this last equation. Note, however, that for $q=1 $, since the mean-value
      of $ u$ is conserved, the smoothness of the flow-map associated with $ \tilde{u} $ ensures the smoothness
      of the flow-map associated with the original equation on hyperplans of functions with a given mean-value.}
       above.
For the cubic Schr\"odinger equation on $ \R $ or $ \T $, this approach shows that the flow-map (if it exists)
cannot be of class $ C^3 $ as soon as $ s<0 $, i.e. $ s_c^\infty=0 $.
A second step was initiated by  Kenig-Ponce and Vega (\cite{KPV4})  who  discovered that more qualitative
ill-posedness phenomena can occur
 above $ s_c $  for dispersive equations
 on the line.  They observed that the lack of uniform continuity on bounded set of the flow-map associated with some
  canonical focusing dispersive equations, including the cubic focusing Schr\"odinger equation, can be proven
   below $ s^\infty_c $ by using a two parameter family of exact solutions.
Using the same idea, Burq-Gerard-Tzvetkov (cf. \cite{BGT1})
noticed that the one-dimensional periodic
 cubic  Schr\"odinger equation (\ref{schro3})
 cannot be uniformly continuous on bounded set  below $ L^2 $.
  The proof  is based on the fact that the solution emanating from the initial data $\alpha e^{inx}$
   is explicitly  given by
 \begin{equation}
 u(t,x)=\alpha \exp\Bigl( -it (n^2-\gamma \alpha^2)\Bigr)\, e^{inx} \; .\label{schrok}
 \end{equation}
  In the same paper, they proved the same type of
   results on the sphere by constructing approximate solutions.

In \cite{CCT1}, Christ-Colliander and Tao also exploited   the construction of approximate solutions in order
 to extend
 the result of \cite{KPV4} in the case of defocusing equations as well as in the periodic case.
Finally, very recently, they proved in \cite{CCT2} an ill-posedness result for the one dimensional periodic cubic Schr\"odinger equation   below $ L^2(\T)  $ in the sense that  the solution map associated with this equation
  (defined on $ L^2(\T) $) cannot be  continuously extended in  $H^s(\T) $ with values in  
    $ C([0,T];H^s(\T)) $ as soon as $ s<0 $.
  The main point is the construction of  approximate solutions corresponding to
   initial data that are supported on two different modes : the mode $ 0 $ plus a mode going to infinity. \\
 In this work we  give a  proof of the ill-posedness of (\ref{schro3}) below $ L^2(\T) $ in the strong sense stated in the beginning of the introduction. Our improvement with respect to the result
  of \cite{CCT2} is that we show the discontinuity of the flow-map for any fixed $t \neq 0$ instead of showing the discontinuity of the  the solution map. Moreover our approach is  
   different and relies on the lack of continuity for the $ L^2(\T)$-weak topology
 of the flow-map associated with this equation. It does not exploit the existence of exact or approximate solutions.
  The main idea is to notice that the nonlinear term of the
   cubic Schr\"odinger
  equation    can be split into two terms. One which enjoys
   a good behavior with respect to the weak convergence in the resolution space. The other which is of the form $ g(\int_{\T}| u|^2) u $. Since $ \int_{\T} |u|^2 $ is a constant of the motion for (\ref{schro3}),
    this will lead to the lack of weak continuity in $ L^2(\T) $ for the flow-map.
     Note that this approach does not work for  the equation  posed on the line. Indeed, it was recently proved in \cite{OL} that  the
  flow-map associated with the cubic Schr\"odinger equation  is continuous in $L^2(\R) $ endowed with the weak topology,
   even if $ s^\infty_c=0 $. It is also interesting to notice that our result uses in  a crucial  way the
    well-posednes
   of the equation in $ L^2(\T) $ which is the critical space for $ C^\infty $-well-posedness.
    Let us note that
   the use of the well-posedness result  in the critical space in order to prove ill-posedness appears also
   in \cite{BT}.  Finally we would like to  mention that we strongly believe that the approach developed
   in this paper can also  lead to the ill-posedness of the periodic Benjamin-Ono equation below $ L^2(\T) $.

\subsection{Main results}
Our ill-posedness result is a straightforward corollary of the following theorem.
\begin{thm} \label{theoshro3}
Let $ \{u_{0,n}\} $ be a sequence of $ L^2(\T) $ converging weakly,  but not
 strongly,  to $ u_0 $ in $ L^2(\T) $ and let $ \{u_n\} $ be the sequence of associated solutions
 of the cubic Schr\"odinger equation (\ref{schro3}).
  For any adherence value $ \alpha^2\in \R_+ $ of $\int_{\T} |u_{0,n}|^2 $
  (at least one such value must be different from   $\int_{\T} | u_{0,n}|^2 $ ) and any increasing sequence of integers
   $ \{n_k\} $ such that $ \int_{\T}|u_{0,n_k} |^2 $ converges to $ \alpha^2 $, the sequence
   $ u_{n_k} (t) $ converges weakly
   in $ L^2(\T) $ to
  $$
     v_\alpha(t) =\exp\Bigl(\frac{i\gamma}{\pi} (\alpha^2-\int_{\T} |u_0|^2) t \Bigr) \, u(t)
     $$
      for all $ t\in \R $, where $ u\in C(\R; L^2(\T))$
   is the  (unique) solution to (\ref{schro3}) emanating from $ u_0 $.
\end{thm}
\begin{rem}
As we already mentioned, in \cite{BGT1} Burg-Gerard-Tzvetkov studied sequences of initial data that are supported on one single
 mode going to infinity whereas in \cite{CCT1}  Christ-Colliander-Tao studied a  sequence of initial data  that are supported on two different modes :
  the  mode $ 0 $ plus a mode going to infinity. These sequences are bounded in $ H^s(\T) $, $s<0$,  but not in $ L^2(\T) $.
  Therefore, we cannot re-examine them  with the help of Theorem \ref{theoshro3}. However,
   we  can
   observe what happens with sequences of initial data having the same support properties. \\
   Sequences of the type
  $ u_{0,n}= \alpha_n e^{i n x} $, with $  \{\alpha_n\}  $ bounded,  converge weakly to $ 0 $ in $ L^2(\T) $
   and Theorem \ref{theoshro3} ensures that the corresponding sequences of emanating solutions tend
    weakly in $ L^2(\T) $ to $0 $ for any $ t\in \R$. Of course this implies the strong convergence to $ 0 $
      in any $ H^s $ for
     $ s<0 $. It is thus clear that this type of counter-examples
    cannot disprove the continuity of the flow-map. \\
Observe now what happens with  a sequence of the type $  u_{0,n}=\beta_1 +\beta_2 e^{inx} $.
 This sequence tends weakly to $ u_0\equiv \beta_1 $ in $ L^2(\T) $ and for any $ n\in \N^*$ it holds $ \|u_{0,n}\|_{L^2}^2=
 2\pi(|\beta_1|^2+|\beta_2|^2) $. The solution of
(\ref{schro3}) emanating from $ u_0=\beta_1 $ is explicitly given by $u(t)=\beta_1
 \, e^{it\gamma |\beta_1|^2 t} $. On the other hand, Theorem \ref{theoshro3} ensures that the
  sequence of emanating solutions
  $ \{u_n\} $ tends weakly in $ L^2(\T) $ and thus strongly in $
  H^s(\T)$, $ s<0 $,
 for any fixed $t\in \R $ toward $ v(t)=\beta_1\, e^{2i\gamma(|\beta_1|^2+ |\beta_2|^2) t} $.
 Since
 $$ |v(t)-u(t)|\ge |\beta_1| \Bigl|1-e^{i 2 \gamma |\beta_2|^2 t} \Bigr|\ge |\gamma| |\beta_1| |\beta_2|^2 t\;
 $$
  for $ t>0 $ small
  enough, this ensures that the
   flow-map is not continuous at $ u_0= \beta_1 $ for the  $ H^s(\T) $-topology as soon as $ s<0 $.
 \end{rem}
 \begin{rem}
We strongly believe  that Theorem \ref{theoshro3} is interesting by its own and not only as a tool to prove ill-posedness.
 For instance, in a forthcoming paper we will use a variant of this theorem to prove the existence of a
 compact  global attractor in $ L^2(\T) $  for the weakly
 damped Schr\"odinger
 equation
 $$
   u_t +i u_{xx}+\varepsilon u +i |u|^2 u=f
 $$
 where $ \varepsilon>0 $ is the damping parameter and where the
external forcing $f$, that is independent of $t$, belongs to
$L^2(\T)$. See \cite{OL} for the case on the line.
\end{rem}

As a consequence of Theorem \ref{theoshro3} we get that  the flow-map of (\ref{schro3}) is not continuous for the weak $ L^2(\T) $-topology. The ill-posedness,
 in the strong sense stated in the introduction, of (\ref{schro3})
 in $ H^s(\T) $, $s<0$, follows directly (see Remark \ref{toto} below).
   \begin{coro} \label{corol}
   For any $ t\neq 0 $, the flow-map associated with the periodic  cubic Schr\"odinger equation (\ref{schro3})
 is discontinuous as a map from $L^2(\T) $, equipped with its weak topology,  into the space of distributions
 $ (C^\infty(\T))^* $ at any point $u_0\in L^2(\T) $ different from the origin.
\end{coro}
\begin{rem} \label{toto}
Let $ B $ be any topological  space  continuously embedded in the space of distributions $ (C^\infty(\T))^* $ and such that $ L^2(\T) $ is compactly embedded in $ B $. Note that this is obviously the case of $ H^s(\T) $ as soon as $ s<0 $.
 Corollary   \ref{corol}   proves the  discontinuity of the map $ u_0\mapsto u(t) $ from $ B $ into $ B $ at any $u_0\in L^2(\T) $ different from the origin  and for
  any fixed $ t\neq 0 $ .
\end{rem}
{\it Proof of Corollary \ref{corol}. } Let $ u_0\in L^2(\T) $ different from $ 0 $ and let $ \{\phi_n\}\subset L^2(\T) $
 be a sequence such that $ \phi_n \rightharpoonup 0 $ in $ L^2(\T) $ 
 and $ \| \phi_n\|_{L^2}^2 \to 2\pi \theta^2 $, $ \theta\in \R^* $, 
   as $ n$ goes to infinity (one can take for instance $ \phi_n=\theta e^{i n x} $). Setting $ u_{0,n}=u_0+
   \phi_n $, we get that $ u_{0,n} \rightharpoonup u_0 $ in $ L^2(\T) $ and
    $ \|u_{0,n} \|_{L^2}^2 \to \|u_0 \|_{L^2}^2+2\pi \theta^2 $ as $
 n\to \infty $. On account of Theorem \ref{theoshro3}, the emanating solutions $ u_n $
  tend weakly in $ L^2(\T) $ for any fixed $ t\in\R$ to
 $ v(t,x)=e^{i 2\gamma \theta^2t} u(t,x) $. Observing  that
 \begin{equation}
 |v(t,x)-u(t,x)|=|1-e^{i 2\gamma \theta^2 t}||u(t,x)|\; \label{titi}
 \end{equation}
 and
 that, $ u_0 \neq 0 $ ensures that  $ u(t)\neq 0 $ for any $ t\in\R   $, we infer that
 $$
  v(t)\neq u(t) \mbox{ in } (C^\infty(\T))^*, \; \forall t\not \in \{\frac{k\pi}{\gamma \theta^2}, \, k\in \Z\} \, ,
 $$
Fixing $ t\neq 0 $ and choosing $ \theta $ such that $ t\not \in \{\frac{k\pi}{\gamma \theta^2}, \, k\in \Z\}$,
  the discontinuity   at $ u_0 $ of the map
   $ u_0 \mapsto u(t) $, from $L^2(\T) $ equipped with its   weak  topology 
   into $ (C^\infty(\T))^* $, follows. \vspace*{2mm}

\section{Proof of Theorem \ref{theoshro3}}
Let us first introduce some notations and function spaces we will work with.
 For  a $ 2\pi $-periodic function $ \varphi$, we define its space Fourier transform  by
$$
\hat{\varphi}(k):=\frac{1}{2\pi} \int_{\T } e^{-i k x} \,
\varphi(x) \, dx , \quad \forall k \in \Z
$$
  and  denote by $ V(\cdot) $ the free group associated with the linearized Schr\"odinger equation,
$$
\widehat{V(t)\varphi}(k):=e^{-i k^2t} \,
\hat{\varphi}(k) , \quad k\in \Z \, .
$$
The Sobolev spaces $ H^s(\T) $ for $ 2\pi$-periodic functions are defined as usually and endowed with
$$
\|\varphi\|_{H^s(\T)}:=\|\langle k \rangle^{s}
\widehat{\varphi}(k) \|_{l^2(\Z)} =\|J^s_x \varphi
\|_{L^2(\T)} \, ,
$$
where $ \langle \cdot \rangle := (1+|\cdot|^2)^{1/2} $ and $
\widehat{J^s_x \varphi}(k):=\langle k \rangle^{s}
\widehat{\varphi}(k)
$. \\
For a function $ u(t,x) $ on $ \T^2 $, we define its space-time Fourier transform by
$$
\hat{u}(q,k):={\mathcal F}_{t,x}(u)(q,k):=\frac{1}{(2\pi)^2}\int_{\T^2}
 e^{-i (q t+ k x)} \, u(t,x) ,
 \quad \forall (q,k) \in  \Z^2 \quad
$$
 and define the  Bourgain spaces $ X^{b,s} $ and ${\tilde X}^{b,s} $ of functions on $ \T^2 $
  endowed with the norm
\begin{equation}
\| u \|_{X^{b,s} } : =
 \| \langle q+k^2\rangle^{b}  \langle k\rangle^s
  \hat{u}\|_{l^2(\Z^2)} =
  \| \langle q\rangle^{b}  \langle k \rangle^s
  {\mathcal F}_{t,x}(V(-t) u ) \|_{l^2(\Z^2)}
\end{equation}
and
\begin{equation}
\| u \|_{{\tilde X}^{b,s} } : =
 \| \langle q-k^2\rangle^{b}  \langle k\rangle^s
  \hat{u}\|_{l^2(\Z^2)} =
  \| \langle q \rangle^{b}  \langle k \rangle^s
  {\mathcal F}_{t,x}(V(t) u ) \|_{l^2(\Z^2)}  \; .
\end{equation}
For a function $ u(t,x) $ on $ \R\times\T $, we define its space-time Fourier transform by
$$
\hat{u}(\tau,\xi):={\mathcal
F}_{t,x}(u)(\tau,\xi):=\frac{1}{2\pi}\int_{\R\times\T} e^{-i (\tau
t+ k x)} \, u(t,x)  dt \, dx , \quad \forall (\tau,k) \in
\R\times\Z \quad .
$$
 and define the  Bourgain spaces $ X^{b,s}_{\R} $ of  functions on $\R\times \T $
  endowed with the norm
$$
\| u \|_{X^{b,s}_{\R} } : =
 \| \langle \tau+k^2\rangle^{b}  \langle k \rangle^s
  \hat{u}\|_{L^2(\R; l^2(\Z))}=
  \| \langle \tau\rangle^{b}  \langle k \rangle^s
  {\mathcal F}_{t,x}(V(-t) u ) \|_{L^2(\R; l^2(\Z))}\; .
$$
Finally,  for $ T>0 $ we define the restriction in time   spaces $ X^{b,s}_T $ of functions
 on $ ]-T,T[ \times \T $
 endowed with the norm
$$
\| u \|_{X_{T}^{b,s}}: =\inf_{v\in X^{b,s}} \{ \|
v\|_{X^{b,s}} , \, v(\cdot)\equiv u(\cdot) \hbox{ on } ]-T,T[
\, \}\;. \vspace*{2mm}
$$
whenever $ 0<T\le 1 $ and
$$
\| u \|_{X_{T}^{b,s}}: =\inf_{v\in X^{b,s}_{\R}} \{ \|
v\|_{X^{b,s}_{\R}} , \, v(\cdot)\equiv u(\cdot) \hbox{ on } ]-T,T[
\, \}\;. \vspace*{2mm}
$$
whenever $ T> 1 $. \vspace{2mm}

 As indicated in the introduction we will use in a crucial way the well-posedness theorem of (\ref{schro3}) in $ L^2(\T) $. So let us recall this theorem proven by Bourgain in \cite{Bo1}.
\begin{thm} \label{Boschro3}
Let $ s\ge 0$. For any $ u_0\in H^s(\T) $ and any $T>0 $, there exists a unique solution
$$
u\in   L^{4}(]-T,T[\times \T)
$$
of (\ref{schro3}). Moreover $ u \in C([-T,T];H^s(\T))\cap X^{1/2,s}_T $ and the map data to solution $ u_0\mapsto u $ is real analytic from $ H^s(\T) $ to $ C([-T,T];H^s(\T)) $.
\end{thm}
This theorem principally use the linear estimates in Bourgain'spaces for the free evolution and the retarded Duhamel operator
\begin{equation}
\| U(t) \varphi\|_{X^{1/2,s}_T} \le C(T) \| \varphi \|_{H^s} \label{freegroup}
\end{equation}
and for any $ 0<\varepsilon <\!< 1 $,
\begin{equation}
\|\int_0^t U(t-t') f(t') \, dt' \|_{X^{1/2,s}_T} \le C(T,\varepsilon) \|f\|_{X^{-1/2+\varepsilon,s}_T}\, ,
\label{duhameloperator}
\end{equation}
 as well as the following periodic estimate to treat the nonlinear term (see \cite{Bo1}):
\begin{equation}
\|v\|_{L^4(\T^2)} \lesssim \|v\|_{X^{3/8,0}}, \quad \forall v\in X^{3/8,0}\, . \label{L4}
\end{equation}
It is worth noticing  that this ensures that for $ 0<T<1 $ it holds
\begin{equation}
\|U(t) \varphi\|_{L^4(]-T,T[\times\T)} \lesssim T^{1/8} \|\varphi \|_{L^2(\T)}, \quad \forall \varphi\in L^2(\T) \, . \label{L41}
\end{equation}
The fact that the time of existence in Theorem \ref{Boschro3} can
be chosen arbitrarly large follows from the conservation of the $
L^2$-norm of the solution. Note finally that excepting the fact
that the solution belongs to $ X^{1/2,s}_T $,
 the remaining of the theorem can be proven for small data 
 by using only the following Zygmund's type  estimate (see again \cite{Bo1}),
 $$
 \| U(t) \varphi\|_{L^4(\T^2)} \lesssim \|\varphi\|_{L^2(\T)} \;,
 $$
 and $ TT^* $ classical arguments.

 As we indicated in the introduction, we are going to split the nonlinear term $ |u|^2 u $
  into different parts.
So let us denote by $ g(\cdot) $ the trilinear operator
\begin{equation}
g(u,v,w):=\bar{u} v w \; \label{trilinear}
\end{equation}
and decompose it in the following way :
\begin{eqnarray*}
g(u,v,w)&=& \sum_{k_1,k_2,k_3\in \Z} \widehat{\bar{u}}(k_1) \widehat{v}(k_2)
\widehat{w}(k_3) e^{i(k_1+k_2+k_3)x} \\
&= &  \sum_{k_1,k_3\in \Z} \widehat{\bar{u}}(k_1) \widehat{v}(-k_1)
\widehat{w}(k_3) e^{i k_3 x} \\
& &+  \sum_{k_1,k_2\in \Z} \widehat{\bar{u}}(k_1) \widehat{v}(k_2)
\widehat{w}(-k_1) e^{i k_2 x} \\
& &-  \sum_{k\in \Z} \widehat{\bar{u}}(k) \widehat{v}(-k)
\widehat{w}(-k) e^{-i k x} \; .\\
&& + \sum_{k_1,k_2,k_3\in \Z\atop (k_1+k_2)( k_1+k_3)\neq 0}
\widehat{\bar{u}}(k_1) \widehat{v}(k_2) \widehat{w}(k_3)
e^{i(k_1+k_2+k_3)x} \; .
\end{eqnarray*}
Note that the three first terms of the above right-hand side correspond to the resonant part of
 $ |u|^2 u $ under the Schr\"odinger flow (see (\ref{resonant}) below).
 Using that $ \widehat{\bar{u}}(k)=\bar{\widehat{u}}(-k)$, it is easy to check for  that, taking $u=v=w$,  the sum of the two first terms gives
  $$
  2\sum_{k_1,k_2\in \Z}| \widehat{u}(k_1)|^2 \widehat{u}(k_2)
 e^{i k_2 x} = \frac{1}{\pi} \|u\|^2_{L^2} u \; .
 $$
 On the other hand, the third term
  is a good term with respect to the weak convergence since its mode $ k$  contains only modes
  $  k  $ or $-k $  of $ u $, $ v $ and $ w$.
We thus rewrite $ g(u):=g(u,u,u) $ as
\begin{equation}
g(u)
   =  \frac{1}{\pi} \|u\|^2_{L^2} u+ \Lambda_1(u,u,u)+\Lambda_2(u,u,u) \;, \label{hihi}
\end{equation}
where
\begin{equation}
\Lambda_1(u,v,w):=\sum_{k_1,k_2,k_3\in \Z\atop (k_1+k_2)(
k_1+k_3)\neq 0} \widehat{\bar{u}}(k_1) \widehat{v}(k_2)
\widehat{w}(k_3) e^{i(k_1+k_2+k_3)x}
\end{equation}
and
\begin{equation}
\Lambda_2(u,v,w):=-  \sum_{k\in \Z} \widehat{\bar{u}}(k) \widehat{v}(-k)
\widehat{w}(-k) e^{-i k x} \; .
\end{equation}
Observe that the conjugaison is an isometry from  $ L^4(\T^2) $
into itself and that,  since  $
\widehat{\overline{v_1}}(q,k)=\overline{\widehat{v_1}(-q,-k)} $
for all $ (q,k)\in\Z^2 $, it is also  an isometry from $
{X^{3/8,0}}$ into ${\tilde X}^{3/8, 0}$. It follows that
(\ref{L4}) still holds  with $X^{3/8,0}$ replaced by ${\tilde
X}^{3/8,0}$ and  that
  \begin{equation} \label{conju}
\|{\mathcal F}^{-1}( |\widehat{\overline{v_1}}|)\|_{L^4(\T^2)}
\lesssim \|\overline{v_1}\|_{{\tilde X}^{3/8,0}} =
\|v_1\|_{{X^{3/8,0}}}\; .
  \end{equation}
Therefore, using suitable extensions of $ u$, $v $ and $ w$ , (\ref{L4}) and (\ref{conju}),
  it is easy to check that for $ i=1,2$ and $ 0<T\le 1  $,
\begin{equation}
\|\Lambda_i(u,v,w)\|_{X^{-7/16,0}_T}\lesssim  \|u \|_{X^{1/2,0}_T}
 \|v \|_{X^{1/2,0}_T}  \|w \|_{X^{1/2,0}_T}
\end{equation}
and thus $ \Lambda_1 $ and $\Lambda_2  $ are continuous operators from $(X^{1/2,0}_T)^3 $ to $ X^{-7/16,0}_T$. Let us now show that they are also continuous  on the same spaces but equipped with their respective weak topology.
\begin{lem}\label{lem1}
The operator $ \Lambda_2 $ is continuous from $(X^{1/2,0}_1)^3 $ to $ X^{-7/16,0}_1$
  equipped with their respective weak topology.
\end{lem}
{\it Proof .} Since $ X^{3/8, -1/3}_1 $ is compactly embedded in $
X^{1/2,0}_1 $,  it suffices to prove that $ \Lambda_2 $ is
bounded from  $ (X^{3/8, -1/3}_1)^3 $ into $ X^{-7/16,-1}_1 $.
But this is straightforward. Indeed, taking extensions $ v_i\in
X^{3/8,-1/3} $ of $ u_i \in X^{3/8,-1/3}_1 $  such that $
\|v_i\|_{X^{3/8,-1/3}}
 \le 2 \|u_i\|_{X^{3/8,-1/3}_1} $,
 it holds
\begin{eqnarray*}
\|\Lambda_2(u_1,u_2,u_3)\|_{X^{-7/16,-1}_1} & \lesssim &
\|\Lambda_2(v_1,v_2,v_3)\|_{X^{-7/16,-1}}  \\
 & \lesssim & \sup_{\| w\|_{X^{7/16,1}}=1} \Bigl| \Bigl(w,\Lambda_2(v_1,v_2,v_3)\Bigr)_{L^2(\T^2)} \Bigr| \; .
\end{eqnarray*}
 Setting  $ q=q_1 +q_2+q_3 $, (\ref{L4}) and (\ref{conju}) then  ensure that
\begin{eqnarray*}
\Bigl| \Bigr(w,\Lambda_2(v_1,v_2,v_3)\Bigl)_{L^2(\T^2)} \Bigr| &
\lesssim   & \sum_{(q_1,q_2,q_3,k)\in \Z^4} |\widehat{w}(q,-k)|
 |\widehat{\overline{v_1}}(q_1,k)| |\widehat{v_2}(q_2,-k)|
|\widehat{v_3}(q_3,-k)| \\
&\ \lesssim & \|{\mathcal F}^{-1}(\langle k \rangle^{-1/3} |\widehat{\overline{v_1}}|)\|_{L^4(\T^2)}
 \|{\mathcal F}^{-1}(\langle k \rangle^{-1/3} |\widehat{v_2}|)\|_{L^4(\T^2)} \\
  & & \hspace*{2mm} \|{\mathcal F}^{-1}(\langle k \rangle^{-1/3} |\widehat{v_3}|)\|_{L^4(\T^2)}
  \|{\mathcal F}^{-1}(\langle k \rangle |\widehat{w}|)\|_{L^4(\T^2)}  \\
 & \lesssim & \|w \|_{X^{3/8,1}}\prod_{i=1}^3  \|v_i \|_{X^{3/8,-1/3}}\;.
 \end{eqnarray*}

\begin{lem}\label{lem2}
The  operator $ \Lambda_1 $ is continuous from $(X^{1/2,0}_1)^3 $ to $ X^{-7/16,0}_1$
  equipped with their respective weak topology.
\end{lem}
{\it Proof .}
For the same reasons as above it suffices to prove that $ \Lambda_1 $ is bounded  from
 $ (X^{7/16,-1/48}_1)^3 $ into $X^{-7/16,-1}_1$.
We proceed as in the proof of the preceding lemma by  introducing extensions $ v_i $ of $ u_i$ such that
 $ \|v_i\|_{X^{3/8,-1/3}} \le 2 \|u_i\|_{X^{3/8,-1/3}_1} $.
Setting $ k=k_1+k_2+k_3 $, we divide the region $
A:=\{(k_1,k_2,k_3)\in \Z^3 \,, (k_1+k_2)(k_1+k_3)\neq 0 \}  $ of
$ \Z^3 $ into  two regions  to estimate
$$
I:=\sum_{(q_1,q_2,q_3)\in \Z^3 \atop (k_1,k_2,k_3) \in A}
|\widehat{w}(q,k)|
 |\widehat{\overline{v_1}}(q_1,k_1)| |\widehat{v_2}(q_2,k_2)|
|\widehat{v_3}(q_3,k_3)|
$$
$\bullet   $ The region $ A_1:=\{(k_1,k_2,k_3)\in A \, , |k| \ge
\frac{1}{4}\max_{i=1,2,3} (|k_i|)\}$. In this region the result
for
 $ \Lambda_1 $  follows exactly the same lines as  for   $ \Lambda_2 $.\\
$\bullet $ The region $ A_2:=\{(k_1,k_2,k_3)\in A \, , |k| <
\frac{1}{4}\max_{i=1,2,3} (|k_i|)\}$.  In this region we first
notice that  $ \max(|k_1+k_2|, |k_1+k_3|)\ge
 \frac{1}{4} \max(|k_i|) $ since otherwise we would have $  \max(|k_1+k_2|, |k_1+k_3|)<
 \frac{1}{4} \max(|k_i|) $ which is clearly incompatible with  $|k|< \frac{1}{4} \max_{i=1,2,3} (|k_i|) $. \\
  Setting $ \sigma=q+k^2, \, \tilde{\sigma}_1=q_1-k_1^2, \, \sigma_2=q_2+k_2^2 $
  and $ \sigma_3=q_3+k_3^2 $ with $ k=k_1+k_2+k_3 $ and $ q=q_1+q_2+q_3 $, we get
  the well-known  resonant relation for the cubic Schr\"odinger equation :
  \begin{equation}
   \sigma-\tilde{\sigma}_1-\sigma_2-\sigma_3 =2(k_1+k_2)(k_1+k_3)  \label{resonant} \;.
   \end{equation}
   This, combining with the fact that $ \min( |k_1+k_2|, |k_1+k_3|)\ge 1 $, ensure that
  $$
  \max(|\sigma|,|\tilde{\sigma}_1|,|\sigma_2|,|\sigma_3|)\gtrsim \max_{i=1,2,3}(|k_i|) \gtrsim \max_{i=1,2,3}\langle k_i\rangle \;.
  $$
It results that
\begin{eqnarray*}
I_{/\Z^3\times A_2 }& \lesssim  & \sum_{(q_1,q_2,q_3)\in \Z^3
\atop (k_1,k_2,k_3)\in A_2} \langle\sigma\rangle^{1/16} |\widehat{w}(q,k)|
 \langle \tilde{\sigma}_1 \rangle^{1/16}\langle k_1\rangle^{-1/48} |\widehat{\overline{v_1}}(q_1,k_1)|  \\
  & & \hspace*{8mm} \prod_{i=2}^3 \langle\sigma_i\rangle^{1/16}  \langle k_i\rangle^{-1/48} |\widehat{v_i}(q_i,k_i)| \\
 & \lesssim & \|\
{\mathcal F}^{-1}( \langle\sigma\rangle^{1/16} |\widehat{w}|)\|_{L^4(\T^2)}
  \|{\mathcal F}^{-1}(\langle\tilde{\sigma}\rangle^{1/16} \langle k \rangle^{-1/48} |\widehat{\overline{v_1}}|)\|_{L^4(\T^2)} \\
& & \hspace*{8mm} \prod_{i=2}^3 \|{\mathcal F}^{-1}(\langle\sigma\rangle^{1/16} \langle k \rangle^{-1/48} |\widehat{v_i}|)\|_{L^4(\T^2)} \\
& \lesssim & \|w \|_{X^{7/16,0}} \prod_{i=1}^3 \|v_i \|_{X^{7/16,-1/48}}
\end{eqnarray*}
where we used (\ref{L4}) and (\ref{conju}) in the last step.\vspace*{2mm} \\

Let $ \{u_{0,n}\}\subset L^2(\T) $ be  a sequence converging weakly but not strongly to $ u_0 $ in $
 L^2(\T) $ and let $\{u_n\} $ be the sequence of emanating solutions. By the Banach-Steinhaus'theorem, 
  $ \{\| u_{0,n}\|_{L^2}\} $ is bounded in $ \R_+ $ and thus admits at
  least one adherence value . Since we assume that $ u_{0,n}  \not \to u_0 $ in $ L^2(\T) $,
  one  adherence value at least  must be different from $ \| u_0\|_{L^2} $.
   Let us denote by  $ \alpha\ge 0 $ such an adherence value of $ \{\| u_{0,n}\|_{L^2}\} $ and let $ \{\| u_{0,n_k}\|_{L^2}\} $  be any subsequence converging towards $ \alpha $.
      From Theorem \ref{Boschro3} we know that the corresponding subsequence of solutions $ \{u_{n_k}\} $ is bounded in $X^{1/2,0}_1 $ and thus, up to the extraction of  a subsequence, converges weakly to some $ v$ in   $X^{1/2,0}_1 $.  Moreover, since the $ L^2 $-norm is conserved for $ u_n $, we infer that
 \begin{eqnarray}
u_n(t) &= & V(t) u_{0,n} - \gamma  \int_0^t V(t-t') \Bigl(\Lambda_1(u_n(t'))+\Lambda_2(u_n(t')\Bigr) \, dt'  \nonumber \\
 &  & -\frac{\gamma}{\pi}\| u_{0,n}\|_{L^2}^2\int_0^t V(t-t') u_n(t')\, dt' \; , \; \forall t \in ]-1,1[ \, .
\end{eqnarray}
From the linear estimates (\ref{freegroup})-(\ref{L4}), Lemmas \ref{lem1} and \ref{lem2}
and the above convergence results, it follows that
\begin{eqnarray}
v(t) &= & V(t) u_{0} - \gamma \int_0^t V(t-t') \Bigl(\Lambda_1(v(t'))+\Lambda_2(v(t')\Bigr) \, dt'  \nonumber \\
 &  & - \frac{\gamma}{\pi}\,  \alpha^2 \int_0^t V(t-t') v(t')\, dt' \;
\end{eqnarray}
and $ v$ is solution of the following Cauchy problem on $ ]-1,1[ $ :
\begin{equation} \label{PDE}
\left\{
\begin{array}{l}
iv_t +v_{xx} +\gamma (\Lambda_1+\Lambda_2)(v)+ \displaystyle \frac{\gamma}{\pi} \alpha^2 v =0 \\
v(0)=u_0
\end{array}
\right. \; .
\end{equation}
Proceeding exactly as for  the cubic Schr\"odinger equation, it is easy to prove that  this Cauchy
 problem is globally well-posed\footnote{Note that the $ L^2 $-norm is  preserved by the flow of (\ref{PDE})} in $ H^s(\T) $, $s\ge 0 $, with a solution belonging  for all $T>0 $ to
$$ C([-T,T];H^s(\T)) \cap L^4(]-T,T[\times \T)\cap X^{1/2,0}_T $$
 with uniqueness in $ L^4(]-T,T[\times \T)$. Therefore, there exists only one possible limit and thus
 the whole sequence $ \{u_{n_k}\} $ converges weakly  to $ v$ in $ X^{1/2,0}_1 $.
Moreover, using the equation satisfied by the $ u_{n} $ and  the
uniform bound in $   L^\infty(]-T,T[; L^2(\T))\cap
L^4(]-T,T[\times \T)$, it is easy to check that for any time-independent $
2\pi $-periodic smooth function $\phi $, the family $ \{t\mapsto
(u_{n_k}(t),\phi)_{L^2} \}
 $ is bounded in $ C([-1,1]) $ and uniformly equi-continuous on $ [-1,1] $. Ascoli's theorem then ensures that $ (u_{n_k}, \phi) $ converges to $(v,\phi)$ on $ [-1,1] $ and thus $ u_{n_k}(t) \rightharpoonup v(t)
  $ in $ L^2(\T) $ for all $ t\in [-1,1] $. By direct iteration  this clearly also holds for all $ t\in \R $.  \\
 Note that, since  the $ L^2 $-norm of the solution is preserved by the flow of
  (\ref{PDE}), $ v$ can be also characterized as the unique solution in $ L^4(]-T,T[\times \T) $ to
\begin{equation}
\left\{
\begin{array}{l}
iv_t +v_{xx} +\gamma |v|^2 v+  \displaystyle \frac{\gamma}{\pi}(\alpha^2- \| u_{0}\|_{L^2}^2) v =0 \\
v(0)=u_0
\end{array}
\right. \; .
\end{equation}
Finally, setting
$$
 \tilde{v}(t,x)=\exp\Bigl( -\frac{i \gamma}{\pi}(\alpha^2- \| u_{0}\|_{L^2}^2)\, t\Bigr)  v(t,x)
 $$
  we notice that $ {\tilde v}\in L^4(]-T,T[\times \T) $ and satisfies
(\ref{schro3}) with $ {\tilde v}(0)=u_0 $. By uniqueness (see Theorem \ref{Boschro3}),  it follows that
${ \tilde v}=u $, where $ u $ is the solution to (\ref{schro3})
 emanating from $ u_0$. This  completes the proof of the theorem.

\noindent
{\bf Acknowledgement:} The author was partially supported by the ANR project "
Etude qualitative des EDP dispersives".


\begin{thebibliography}{99}
\bibitem{A}  M. Ablowitz and Y. Ma, {\em The periodic nonlinear Schr\"odinger equation}, Stud. Appl. Math. 65
 (1981), 113-158.

\bibitem{BT} I. Bejenaru and T. Tao, {\sl  Sharp well-posedness and ill-posedness results for a quadratic non-linear {S}chr\"odinger equation.}  J. Funct. Anal.  233  (2006),  no. 1, 228--259.

\bibitem {Bo1} J.~Bourgain,
{\em  Fourier transform restriction phenomena for certain lattice subsets and application to nonlinear evolution equations
I. The Schr\"odinger equation},  GAFA, 3 (1993), 157-178.

\bibitem {Bo4}  J.~Bourgain,
{\em Periodic Korteweg de Vries equation with measures as initial data}
, Sel. Math. New. Ser. 3 (1993), pp. 115--159.

\bibitem{BGT1} N. Burq, P. G\'erard and N. Tzvetkov, {\sl An instability property of the nonlinear {S}chr\"odinger equation on $S\sp d$},  Math. Res. Lett.  9  (2002),  no. 2-3, 323--335.


\bibitem{CCT1} M. Christ, J. Colliander and T. Tao,
 {\sl Asymptotics, frequency modulation, and low regularity ill-posedness for canonical defocusing equations},
Amer. J. Math. 125 (2003), no. 6, 1235--1293.

\bibitem{CCT2} M. Christ, J. Colliander and T. Tao,
 {\sl Instability of the periodic nonlinear {S}chr\"odinger equation}, preprint arXiv:math/0311048.


\bibitem{OL} O. Goubet and L. Molinet, {\sl Global weak attractor for weakly damped nonlinear
 {S}chr\"odinger equations in $ L^2(\R) $}, submitted.




\bibitem{KPV4} C. E. Kenig, G. Ponce and L. Vega, {\sl  On the ill-posedness of some canonical dispersive equations},  Duke Math. J.  106  (2001),  no. 3, 617--633.

\bibitem{Ts} Y. Tsutsumi, { \sl  $ L^2 $-solutions for nonlinear {S}chr\"odinger equations and nonlinear groups},
   Funk. Ekva. 30  (1987),  115-125.









\end{thebibliography}
\end{document}